\documentclass[12pt]{article}
\usepackage{amssymb,amsfonts}
\usepackage{graphicx}
\usepackage{amscd}
\newtheorem{theorem}{Theorem}

\title{Invariant Lagrange Submanifolds of Dissipative Systems}
\author{A. Agrachev \thanks{SISSA, Trieste and MIAN, Moscow}}
\date{}
\begin{document}
\maketitle

\begin{abstract} We study smooth solutions of modified
Hamilton--Jacobi equations $H(\frac{du}{dq},q)+\alpha u(q)=0,\
q\in M$, on a compact manifold $M$.
\end{abstract}

 Let $M$ be a compact
Riemannian manifold of class $C^k,\ k\ge 2,$ with the Riemannian
structure $(\xi,\eta)\mapsto\langle I_q^{-1}\xi,\eta\rangle,\
\xi,\eta\in T_qM,\ q\in M$, where $I_q:T^*_qM\to T_qM$ is a
self-adjoint linear map such that the quadratic form
$z\mapsto\langle z,I_qz\rangle,\ z\in T_q^*M,$ is positive definite.

Let $V\in C^k(M)$ and $\omega$ be a closed differential 1-form on
$M$ of class $C^k$ such that $\nabla\omega=0$, where
$\nabla\omega$ is the covariant derivative of $\omega$. We
consider the Hamiltonian function $H\in C^k(T^*M)$ defined by the
formula:
$$
H(z)=\frac 12\langle I_q(z+\omega_q),z+\omega_q\rangle+V(q),\quad
z\in T^*_qM.
$$
Let $\vec H$ be the Hamiltonian vector field on $T^*M$ associated
to $H$ and $\ell$ be the ``vertical" Euler vector field of the
vector bundle $T^*M\to M$. In local coordinates, $z=(p,q),\
p,q\in\mathbb R^n,\ T^*_qM=(\mathbb R^n,q)$,
$$
H(p,q)=\frac12(p+\omega(q))^*I_q(p+\omega(q))+V(q),
$$
$\vec H(p,q)=\sum\limits_{i}\left(\frac{\partial H}{\partial
p^i}\frac\partial{\partial q^i}-\frac{\partial H}{\partial
q^i}\frac\partial{\partial p^i}\right),\
\ell(p,q)=\sum\limits_{i}p^i\frac\partial{\partial p^i}. $

We study the dissipative system $\dot z=\vec H(z)-\alpha\ell(z),$
where $\alpha$ is a positive constant.
It is not hard to see that any bounded trajectory of this system
is contained in the set
$$
B_H\stackrel{def}{=}\{z\in T^*M:
H(z-\omega_{\pi(z)})\le\max\limits_{q\in M}H(0_q)\},
$$
where $0_q$ is the origin of $T^*_qM$ and $\pi:T^*M\to M,\ 
\pi(T^*_qM)=q$.

Given $z\in T^*M$, we denote by $\rho(z)$ the maximal eigenvalue
of the symmetric operator
$$
\xi\mapsto\mathfrak R(\xi,I_qz)I_qz+(\nabla^2_qV)\xi,\quad \xi\in
T_qM,
$$
where $\mathfrak R$ is the Riemannian curvature. Finally, we set
$$r=\max\{\rho(z): z\in B_H\}.$$

Let $\Omega^\alpha$ be the set of all absolutely continuous curves
 $\gamma:[0,+\infty)\to M$ such that the integral
$\int\limits_0^{+\infty}e^{-\alpha t}\langle
I^{-1}_{\gamma(t)}\dot\gamma(t),\dot\gamma(t)\rangle\,dt$
converges.
We introduce the {\it discounted action} functional
$$
\mathfrak I_\alpha(\gamma)= \int\limits_0^{+\infty}e^{-\alpha
t}\left(\frac 12\langle
I^{-1}_{\gamma(t)}\dot\gamma(t),\dot\gamma(t)\rangle-V(\gamma(t))+
\langle\omega_{\gamma(t)},\dot\gamma(t)\rangle\right)\,dt,\quad
\gamma\in\Omega_\alpha.
$$

\begin{theorem} Let $u(q)=-\inf\{\mathfrak
I_\alpha(\gamma):\gamma\in\Omega_\alpha,\gamma(0)=q\},\ q\in M$.
If $r\le 0$ or $0<r<\frac{\alpha^2}4$ and $k<\frac
2{1-\sqrt{1-\frac{4r}{\alpha^2}}}$, then:
\begin{itemize}
\item $u\in C^k(M)$ and the map $(H,\alpha)\mapsto u$ is continuous
in the $C^2$-topology.
\item The function $u$ satisfies the
modified Hamilton--Jacobi equation \linebreak $H(du)+\alpha u=0$
and $\{d_qu:q\in M\}\subset T^*M$ is an invariant submanifold of
the system $\dot z=\vec H(z)-\alpha\ell(z)$.
\item There exists a containing 0 neighborhood $\mathcal O$ of $u$ in $C^2(M)$
such that $\forall v_0\in\mathcal O$ the classical solution $v_t$
of the Cauchy problem $\frac{\partial u_t}{\partial
t}+H(du_t)+\alpha u_t=0,\ u_0=v_0$, is defined for all $t\ge 0$
and $\|dv_t-du\|_{C^1}\to 0$ as $t\to +\infty$ with the
exponential convergence rate.
\end{itemize}
\end{theorem}

\noindent{\bf Remark.} Theorem 1 is applied to the Hamiltonians in
$\mathbb R^n\times\mathbb R^n$ of the form $$ H(p,q)=\frac
12|p+a|^2+V(q),$$ where $a\in\mathbb R^n$ is a constant vector and
$V$ is a smooth periodic potential. Then $r$ is the maximum of the
eigenvalues of the matrices $\frac{d^2V}{dq^2},\ q\in\mathbb R^n$.
If $r<\frac{\alpha^2}4$, then the equation $ \frac
12|\frac{du}{dq}+a|^2+V(q)+\alpha u=0$ has a periodic
$C^k$-solution $u$, where $k$ is maximal integer that is strictly
smaller than $\frac 2{1-\sqrt{1-\frac{4r}{\alpha^2}}}$. Moreover,
$\{(\frac{du}{dq},q):q\in\mathbb R^n\}$  is an invariant
submanifold of the system $$ \dot q=p+a, \quad \dot
p=-\frac{dV}{dq}-\alpha p.$$

\medskip

The proof of Theorem~1 can be derived from \cite{ag} and
\cite{HPS}. Indeed, Theorem~1 is an improvement of results of
paper \cite{ag}. The improvement concerns more general
Hamiltonians (nonzero forms $\omega$ are available), better
smoothness of $u$, and stability properties. One can check that
$\omega$ does not affect the canonical connection and the
curvature operators; hence more general Hamiltonians do not
require essential changes in the proof.

The better smoothness and stability follow from
\cite[Th.4.1]{HPS}. Indeed, Prop.~1 in \cite{ag} implies that
$\{d_qu:q\in M\}$ is a normally hyperbolic invariant submanifold
(see \cite{HPS} for the definition) of the flow generated by the
vector field $\vec H(z)-\alpha\ell(z)$. Moreover, this normally
hyperbolic invariant submanifold has zero unstable subbundle and
can be actually called the ``normally stable" invariant
submanifold. Theorem~4.1 in \cite{HPS} contains the estimate for
the smoothness class of the normally hyperbolic invariant
submanifold in terms of the Lyapunov exponents while the
analysis of the proof of Prop.~1 in \cite{ag} gives explicit
estimates for the Lyapunov exponents in terms of $r$ and
$\alpha$.


\begin{thebibliography}{9}

\bibitem{ag} A. Agrachev, \textit{Well-posed infinite horizon variational
problems on a compact manifold}, arXiv:0906.4433, 23pp.

\bibitem{HPS} M. Hirsch, C. Pugh, M. Shub, \textit{Invariant
manifolds}, Lecture Notes in Math., \textbf{583}. Springer Verlag,
1977, 149pp.

\end{thebibliography}
\end{document}